\documentclass{article}

\usepackage{delarray,verbatim,enumerate,a4wide}
\usepackage{graphicx}
\usepackage{xcolor}
\usepackage{amsmath,amsthm,amstext,amsbsy,amssymb,amsfonts,amscd}
\usepackage[all]{xy}
\usepackage[frenchb]{babel}
\usepackage[T1]{fontenc}
\usepackage[utf8]{inputenc}

\usepackage{upgreek}
\usepackage{color}
\definecolor{vert}{rgb}{0.1,0.4,0.2}
\usepackage[colorlinks=true,linkcolor=blue,citecolor=vert]{hyperref}

\usepackage{calligra}
\DeclareFontShape{T1}{calligra}{m}{n}{<->s*[0.95]callig15}{}
\DeclareMathAlphabet{\mathscr}{T1}{calligra}{m}{n}

\newtheorem{Th}{Théorème}[]
\newtheorem{Lem}[Th]{Lemme}
\newtheorem{Prop}[Th]{Proposition}
\newtheorem{Cor}[Th]{Corollaire}

\newtheorem{Sco}[Th]{Scolie}
\newtheorem{Def} [Th]{Définition}

\def\Preuve{\smallskip\noindent {\it Preuve.~}}

\font\teneufm=eufm10
\font\seveneufm=eufm7
\font\fiveeufm=eufm5
\newfam\mathfrakfam
\textfont\mathfrakfam=\teneufm \scriptfont\mathfrakfam=\seveneufm
\scriptscriptfont\mathfrakfam=\fiveeufm
\def\mathfrak{\fam\mathfrakfam}

		\def\QQ{\mathbb Q}	
\def\NN{\mathbb N}	\def\ZZ{\mathbb Z}
\def\F2{\mathbb{F}_2}	\def\Z2{\mathbb{Z}_2}		
\def\Zl{\mathbb{Z}_\ell} 	\def\Ql{\mathbb{Q}_\ell}	\def\Tl{\mathbb{T}_\ell}

 					\def\F{\mathcal  F}
  		\def\R{\mathcal  R}

		\def\x{{\mathfrak x}}

\def\Gal{\operatorname{Gal}}		\def\Rad{\operatorname{Rad}}
		\def\Hom{\operatorname{Hom}}

\makeatletter
\newcommand*\wt[2][0.2ex]{%
        \begingroup
        \mathchoice{\wt@helper{#1}{#2}{\displaystyle}{\textfont}}
                   {\wt@helper{#1}{#2}{\textstyle}{\textfont}}
                   {\wt@helper{#1}{#2}{\scriptstyle}{\scriptfont}}
                   {\wt@helper{#1}{#2}{\scriptscriptstyle}{\scriptscriptfont}}%
        \endgroup
        #2%
}
\newcommand*\wt@helper[4]{%
        \def\currentfont{\the#41}%
        \def\currentskewchar{\char\the\skewchar\currentfont}%
        \setbox\tw@\hbox{\currentfont$#2$\currentskewchar}%
        \dimen@ii\wd\tw@
        \setbox\tw@\hbox{\currentfont$#2${}\currentskewchar}%
        \advance\dimen@ii-\wd\tw@
        \rlap{\raisebox{-#1}{$\m@th#3\kern\dimen@ii\widetilde{\phantom{#2}}$}}%
}
\makeatother

\def\wE{\,\wt[0.2ex]{\!\mathcal E}}

\newcommand{\Bmu}{\mbox{$\raisebox{-0.59ex}
{$l$}\hspace{-0.18em}\mu\hspace{-0.88em}\raisebox{-0.98ex}{\scalebox{2}
{$\color{white}.$}}\hspace{-0.416em}\raisebox{+0.88ex}
 {$\color{white}.$}\hspace{0.46em}$}{}}
\newcolumntype{x}[1]{>{\centering\hspace{0pt}}p{#1}}

\begin{document}

\title{\LARGE\bf  Sur le radical kummérien des $\Zl$-extensions}

\author{ Jean-François {\sc Jaulent} }
\date{}
\maketitle
\bigskip\bigskip

{\small
\noindent{\bf Résumé.}
Sur la base d'un travail antérieur, nous donnons une description nouvelle du radical initial attaché au compositum des $\Zl$-extensions d'un corps de nombres en termes de limites projectives pour la norme dans la $\Zl$-extension cyclotomique de ce corps. Le résultat précis que nous obtenons contient, comme conséquence immédiate, les encadrements de ce radical proposés par Soogil Seo dans une série d'articles.  Nous montrons enfin qu'une description analogue vaut pour le noyau des symboles universels de Tate.}

\smallskip

{\small
\noindent{\bf Abstract.}
On the basis of a previous work, we elaborate a new description of the Kummer radical associated to the first layers of $\Zl$-extensions of a number fields $K$, by using inverse limits for the norm maps in the cyclotomic $\Zl$-extension $K_\infty/K$. Our main result contains, as an obvious consequence, the inclusions provided by Soogil Seo in a set of papers. By the same way we also give in the last section a similar description of the Tate kernel for universal symbols in $K_2(K)$.}
\bigskip\bigskip


\bigskip

\noindent{\bf  Introduction}\medskip

La détermination du radical de Kummer attaché au compositum des $\Zl$-extensions d'un corps de nombres pose des problèmes difficiles, tant théoriques que pratiques: au plan théorique, car ils mettent en jeu de façon non triviale la conjecture de Leopoldt; au plan pratique parce qu'ils nécessitent le plus souvent de conduire des calculs dans les divers étages de la tour cyclotomique, ce qui augmente rapidement le degré des corps sur lesquels on est amené à travailler.\smallskip

De multiples travaux ont ainsi été publiés ces dernières années, qui traitent de ces questions dans le cadre de la théorie d'Iwasawa ou de la théorie du corps de classes (ainsi \cite{Gr1}, \cite{Grt}, \cite {J18}, \cite{J23}, \cite{J24} et \cite{J31}), éventuellement en lien avec la $K$-théorie des corps de nombres  (comme \cite{JM} et \cite{MN2}) ou en proposent des approches plus effectives  (ainsi  \cite{Gr4}, \cite{He} et \cite{Th}) illustrées par de nombreux exemples.\smallskip

Plus récemment Seo, dans une série d'articles \cite{Seo1, Seo2,Seo3} qui méconnaissent malheureusement la totalité de ceux cités plus haut, a mis en avant divers encadrements de ce même radical.
\medskip

L'objet de la présente note est de montrer directement comment un lemme technique très simple, déjà utilisé par Bertrandias et Payan dans leur article fondateur \cite{BP}, permet de déduire de la description kummérienne exposée dans \cite{J18} et \cite{J23} une identité qui interprète le radical initial des $\Zl$-extensions en termes de limites projectives pour les applications normes (Théorème \ref{RadZ} et Scolie) et entraîne de façon immédiate les approximations successives avancées par Seo.\smallskip

En complément, nous montrons comment la même technique permet d'atteindre le noyau des symboles universels de Tate.

\bigskip

\noindent{\bf  Conventions}\medskip


Dans tout ce qui suit $\ell$ désigne un nombre premier arbitraire, $K$ un corps de nombres qui contient les racines $2\ell$-ièmes de l'unité et $K_\infty=\cup_{n\in\NN} K_n$ sa $\Zl$-extension extension cyclotomique.\smallskip

On note $\ell^m=|\mu^{\phantom{lc}}_K|$ le cardinal du $\ell$-groupe des racines de l'unité contenues dans $K$ (avec donc $m>1$ pour $\ell=2$ et $m>0$ sinon), de sorte que le $n$-ième étage de la tour, de degré $[K_n:K]=\ell^n$, est le corps $K_n=K[\Bmu^{\phantom{lc}}_{\ell^{n+m}}]$ engendré sur $K$ par les racines $\ell^{m+n}$-ièmes de l'unité.\smallskip

On désigne enfin par $\gamma$ le générateur topologique du groupe procyclique $\Gamma=\Gal(K_\infty/K)$ défini par la formule $\zeta^\gamma=\zeta^{1+l^m}$, pour tout $\zeta\in\Bmu^{\phantom{lc}}_{\ell^\infty}$; et par $\Gamma^*=\gamma^{*\,\Zl}$ le groupe multiplicatif topologiquement engendré par l'élément $\gamma^*=\gamma/(1+\ell^m)$ de l'algèbre d'Iwasawa $\Lambda=\Zl[[\gamma-1]]$.

\newpage

\newpage

\noindent{\bf 1. Rappels sur la théorie de Kummer cyclotomique}\medskip


Commençons par rappeler la définition du $\ell$-radical universel (cf. \cite{J18},  I.2, ou \cite{J23}, Déf. 1.1):

\begin{Def} Par $\ell$-radical universel d'un corps de nombres $K$, nous entendons le tensorisé par $\QQ_{\ell}/\Zl$ de son groupe multiplicatif, que nous notons:\smallskip

\centerline{${\mathfrak R}_K=(\mathbb Q_{\ell}/\mathbb Z_{\ell})\otimes_{\ZZ} K^\times$.}
\end{Def}

Le groupe ${\mathfrak R}_K$ est, par construction, un $\mathbb Z_{\ell}$-module divisible de torsion, dont chaque élément admet une représentation de la forme $\ell^{-k}\otimes x$ avec $k\in\NN$ et $x\in K^\times$.
La condition $\ell^{-k}\otimes x=1$ signifie alors que l'élément $x$ est le produit $x=\zeta \ y^{\ell^k}$ d'une racine de l'unité et d'une puissance $\ell^k$-ième dans $K^\times$. On prendra garde, en effet, que le produit tensoriel par le groupe divisible $\mathbb Q_\ell/\mathbb Z_\ell$ a pour effet de tuer le sous-groupe de torsion de $K^\times$.\smallskip

Dans le cas qui nous intéresse, la $\Zl$-extension cyclotomique $K_\infty$ étant la réunion des sous-corps $K_n$, le radical universel ${\mathfrak R}_{K_\infty}$ qui lui correspond s'obtient naturellement comme limite inductive des radicaux ${\mathfrak R}_{K_n}=(\QQ_\ell/\ZZ_\ell)\otimes_{\ZZ} K^\times_n$. Plus précisément (cf. \cite{J18},  I.2, ou \cite{J23}, Prop. 1.2), il vient:

\begin{Prop} L'application naturelle de $\,{\mathfrak R}_{K_n}$ dans $\,{\mathfrak R}_{K_\infty}$ est un morphisme
injectif. Et, pour chaque entier naturel $n$, le radical universel $\,{\mathfrak R}_{K_n}$ associé au corps
$K_n$ s'identifie canoniquement au sous-module des points fixes dans $\,{\mathfrak R}_{K_\infty}$ du groupe de Galois
$\Gamma_{\!n}=Gal(K_\infty/K_n)$~; ce qui s'écrit tout simplement:

\centerline{${\mathfrak R}_{K_n}^{\phantom{l}}=\,{\mathfrak R}^{\,\Gamma_{\!n}}_{K_\infty},\quad\forall n\in\mathbb N$.}
\end{Prop}

Le corps $K_\infty$ contenant par construction le groupe $\Bmu^{\phantom{lc}}_{\ell^\infty}$ de toutes les racines de  l'unité d'ordre $\ell$-primaire, le groupe $\,{\mathfrak R}_{K_\infty}$ permet de décrire par radicaux l'ensemble des pro-$\ell$-extensions abéliennnes de $K_\infty$, le radical de Kummer attaché à une telle extension $L_\infty/K_\infty$ étant donné par:\smallskip

\centerline{$\Rad (L_\infty/K_\infty)=\{\ell^{-k}\otimes x_\infty \in {\mathfrak R}_{ K_\infty}\;|\; K_\infty[\sqrt[\ell^k]{x_\infty}\,] \subset L_\infty\}$;}\smallskip

\noindent et son groupe de Galois par:\smallskip

\centerline{$\Gal(L_\infty/K_\infty)=\Hom (\Rad (L_\infty/K_\infty),\Bmu^{\phantom{lc}}_{\ell^\infty})$.}\smallskip

En particulier, les pro-$\ell$-extensions $L_\infty/K_\infty$ qui proviennent par composition avec $K_\infty$ d'une pro-$\ell$-extension abélienne $L/K$ sont celles dont le groupe de Galois est invariant par $\Gamma=\gamma^{\Zl}$, c'est-à-dire celles dont le radical de Kummer est invariant par $\Gamma^*=\gamma^{*\,\Zl}$.
Plus généralement:

\begin{Th}\label{Abel}
Les pro-$\ell$-extensions $L_nK_\infty$ de $K_\infty$ qui proviennent d'une pro-$\ell$-extensions abélienne $L_n$ de $K_n$ sont celles dont le radical de Kummer est invariant par $\Gamma^*_{\!n}=\gamma^{*\,\ell^n \Zl}$. Si donc $K^{ab}_n$ désigne la pro-$\ell$-extension abélienne maximale de $K_n$, il suit:\smallskip

\centerline{$\Rad(K^{ab}_nK_\infty/K_\infty)= {\mathfrak R}_{K_\infty}^{\,\Gamma^*_{\!n}}$.}\smallskip

\noindent Et le radical de Kummer associé au compositum $Z_n$ des $\Zl$-extensions de $K_n$ est ainsi son sous-groupe $\ell$-divisible maximal:

\centerline{$\Rad(Z_nK_\infty/K_\infty)= ({\mathfrak R}_{K_\infty}^{\Gamma^*_{\!n}})_{{}_{\mathrm{div}}}$.}
\end{Th}

\begin{Sco}
Le radical kummérien de la sous-extension maximale de  $K^{ab}_n$ qui peut être décrite radiciellement par des éléments de $K_n$ est le sous-groupe de $\ell^{m+n}$-torsion du radical universel $\,{\mathfrak R}_{K_n}$:\smallskip

\centerline{${\mathfrak R}^{\,\Gamma_{\!n}}_{K_\infty}\cap\;{\mathfrak R}^{\,\Gamma^*_{\!n}}_{K_\infty} = {}_{\ell^{m+n}}{\mathfrak R}_{K_n}^{\phantom{l}}$.}
\end{Sco}

\Preuve Le radical cherché est formé, en effet, des éléments de ${\mathfrak R}_{K_\infty}$ qui sont simultanément dans ${\mathfrak R}^{\,\Gamma_{\!n}}_{K_\infty}={\mathfrak R}_{K_n}^{\phantom{l}}$ et dans ${\mathfrak R}^{\,\Gamma^*_{\!n}}_{K_\infty}$, c'est-à-dire ceux fixés à la fois par $\gamma^{\ell^n}$ et par $\gamma^{*\ell^n}= \gamma^{\ell^n}/(1+\ell^m)^{\ell^n}$, donc par $(1+\ell^m)^{\ell^n}$ et  finalement par $(1+\ell)^{m+n}$, puisque le quotient $(1+\ell^m)^{\ell^n} / (1+\ell^{m+n})$ est inversible dans $\Zl$. C'est essentiellement le corollaire 1.3 de \cite{J31}. \medskip

On retrouve le fait banal que les $\ell$-extensions abéliennes d'un corps de nombres $K$ qui peuvent être décrite radiciellement à partir des seuls éléments de $K$ sont exactement celles dont le groupe de Galois est d'exposant $\ell^m$, dès lors que le corps $K$ contient les racines $\ell^m$-ièmes de l'unité mais non  les racines $\ell^{m+1}$-ièmes

\newpage

\noindent{\bf 2. Description radicielle du compositum des $\Zl$-extensions}\medskip


Pour établir le résultat principal de cette note, partons d'un lemme purement algébrique, dont le premier item est également le point de départ de Bertrandias et Payan (cf. \cite{BP}, Lemme 1):

\begin{Lem}\label{LemmeBP} Soient $\gamma$ une indéterminée, $\Zl[\gamma]$ l'algèbre des polynômes 
à coefficients dans $\Zl$ et $\gamma^*=\gamma/(1+\ell^m)$ avec $m>1$ pour $\ell=2$ et $m>0$ sinon. Pour tout $n\in\NN$, on a  les congruences:\smallskip

\begin{itemize}
\item[(i)] $\nu_n = 1+\gamma+\gamma^2+\cdots+\gamma^{\ell^n-1} \equiv u_n\ell^n \quad  [\!\!\mod (\gamma^*-1)\,]$, avec $u_n\in 1+\ell\ZZ$;\smallskip

\item[(ii)] $\omega_n=\gamma^{\ell^n}-1 \equiv  v_n\ell^{m+n} \quad  [\!\!\mod (\gamma^*-1)\,]$, avec $v_n \in 1+\ell\ZZ$;\smallskip

\item[(iii)] $\omega_{n+1}/\omega_n = 1+\gamma^{\ell^n}+ \cdots + \gamma^{\ell^n(\ell-1)} \equiv \ell v'_n  \quad  [\!\!\mod (\gamma^*-1)\,]$, avec $v'_n \in 1+\ell\ZZ$.

\end{itemize}
\end{Lem}

\Preuve Raisonnons modulo $(\gamma^*-1)$. De la congruence $\gamma \equiv 1+\ell^m$, nous tirons successivement:
\begin{itemize}
\item[(i)] $\nu_n \equiv 1 + (1+\ell^m) + (1+\ell^m)^2 + \cdots + (1+\ell^m)^{\ell^n-1} \equiv \frac{(1+\ell^m)^{\ell^n}-1}{(1+\ell^m)-1} \equiv u_n \ell^n$, où $u_n$ est dans $1+\ell\ZZ$ (donc, en particulier, inversible dans $\Zl)$;
\item[(ii)] $\omega_n=\gamma^{\ell^n}-1 \equiv (1+\ell^m)^{\ell^n}-1 \equiv v_n\ell^{m+n}$, pour un $v_n$ dans $1+\ell\ZZ$;\smallskip

\item[(iii)] $\omega_{n+1}/\omega_n = 1+\gamma^{\ell^n}+ \cdots + \gamma^{\ell^n(\ell-1)} \equiv 1 + (1+v_n\ell^{m+n}) + \cdots + (1+v_n\ell^{m+n})^{\ell-1}$, d'où: $\omega_{n+1}/\omega_n \equiv \frac{(1+v_n\ell^{m+n})^\ell-1}{(1+v_n\ell^{m+n})-1} \equiv v'_n\ell$, pour un $v'_n$ dans $1+\ell\ZZ$, comme annoncé.
\end{itemize}\medskip

Intéressons-nous maintenant au radical initial du compositum $Z$ des $\Zl$-extensions du corps $K$, c'est-à-dire, puisque $\mu_K^{\phantom{lc}}$ est d'ordre  $\ell^m$, aux éléments de $\ell^m$-torsion du radical $({\mathfrak  R}_{K_\infty}^{\,\Gamma^*})_{{}_{\mathrm{div}}}$ associé par le Théorème \ref{Abel} à la pro-$\ell$-extension abélienne $ZK_\infty/K_\infty$.\par

Soit donc $\x\in{}_{\ell^m}^{\phantom{l}}\!{\mathfrak  R}_K^{\,\Gamma^*}
\cap\,({\mathfrak  R}_{K_\infty}^{\,\Gamma^*})_{{}_{\mathrm{div}}}$ un tel élément; 
c'est la puissance $\ell$-ième $\x=\x_1^\ell$ d'un élément $\x_1$ de ${}_{\ell^{m+1}}{\mathfrak  R}_{K_1}^{\,\Gamma^*}$, qui est lui-même la puissance $\ell$-ième $\x_1=\x_2^\ell$ d'un élément $\x_2$ de ${}_{\ell^{m+2}}{\mathfrak  R}_{K_2}^{\,\Gamma^*}$ et ainsi de suite. En d'autres termes, le radical cherché est l'ensemble des termes initiaux des suites cohérentes $(\x_n)_{n\in\NN}$ ainsi construites, i.e. des éléments de la limite projective $\varprojlim \,{\mathfrak  R}_{K_n}^{\,\Gamma^*}$ pour l'exponentiation.\par

Cela étant, tous les $\x_n$ sont dans ${\mathfrak  R}_{K_\infty}^{\,\Gamma^*}$, qui est un $\Zl[\gamma]$-module annulé par $\gamma^*-1$. Et d'après le dernier item du Lemme, comme opérateur sur un tel module, l'exponentiation par $\ell$, à un inversible près, agit comme la norme $\omega_{n+1}/\omega_n$. Nous pouvons donc modifier de proche en proche les $\x_n$ pour remplacer les identités de divisibilité $\x_n=\x_{n+1}^\ell$ par les identités normiques  $\x_n=\x_{n+1}^{\omega_{n+1}/\omega_n}$. Ainsi:

\begin{Th}\label{RadZ}
 Le radical initial attaché au compositum $Z$ des $\Zl$-extensions de $K$, i.e. le sous-groupe du radical universel $\,{\mathfrak R}_{K_\infty}$ attaché au compositum des $\ell$-extensions de $K_\infty$ qui proviennent par composition avec $K_\infty$ d'une extension cyclique d'exposant $\ell^m$ de $K$ qui est $\Zl$-plongeable\smallskip

\centerline{$_{\ell^m\!}{\mathfrak Z}_K = \{\ell^{-m}\otimes x \in {\mathfrak R}_K \; | \; K_\infty [\sqrt[\ell^m]{x}] \subset Z\}$}\smallskip 

\noindent est formé des termes initiaux $\x_0$ de la limite projective pour les applications normes $\varprojlim \,{\mathfrak  R}_{K_n}^{\,\Gamma^*}$.
\end{Th}

Pour interpréter les radicaux $\,{\mathfrak  R}_{K_n}^{\,\Gamma^*}$ en termes de quotients, introduisons les $\ell$-adifiés\smallskip

\centerline{$\R_{K_n}=\Zl \otimes_{\ZZ} K_n^\times$}\smallskip

\noindent des groupes multiplicatifs des corps $K_n$, autrement dit les groupes d'idèles principaux attachés aux $K_n$ par la théorie $\ell$-adique du corps de classes (cf. e.g. \cite{J18} ou \cite{J31}). Et notons $V_{K_n}=\Ql\otimes_{\Zl}\R_{K_n}=\Ql\otimes_{\ZZ}K^\times_n$ les $\Ql$-espaces associés. Formons enfin pour $n$ fixé la suite exacte courte:\smallskip

\centerline{$ 1 \rightarrow \R^{\phantom{l}}_{K_n}/\mu_{K_n}^{\phantom{lc}} \rightarrow V^{\phantom{l}}_{K_n} \rightarrow {\mathfrak R}^{\phantom{l}}_{K_n} \rightarrow 1$;}\smallskip

\noindent et considérons le diagramme du serpent associé à l'opérateur $\gamma^*-1$. Observons que les $V_{K_n}$ sont les réunions des $V_{K_n}^S=\Ql\otimes_{\ZZ} E^{S}_{K_n}$ obtenus à partir des $S$-unités lorsque $S$ parcourt les ensembles finis de places de $K$ et que les $V_{K_n}^S$ sont des $\Ql$-espaces vectoriels de dimension finie, stables par $\gamma$ et fixés par $\gamma^{\ell^n}$. Ainsi $\gamma^*-1$ opère bijectivement sur les $V_{K_n}^S$ donc sur $V_{K_n}$; et il suit:\smallskip

\centerline{${\mathfrak  R}_{K_n}^{\,\Gamma^*} \simeq \R^{\phantom{l}}_{K_n}/\R_{K_n}^{\gamma^*-1}\mu_{K_n}^{\phantom{lc}}$.}

\begin{Sco}
Le radical initial attaché au compositum des $\Zl$-extensions de $K$ s'identifie au terme initial de la limite projective pour les applications normes des quotients des groupes d'idèles principaux invariants par $\gamma^*-1$ attachés aux étages finis de la tour cyclotomique $K_\infty/K$:
\smallskip

\centerline{$_{\ell^m\!}{\mathfrak Z}^{\phantom{l}}_K = \pi_0 \left( \varprojlim \;(\R^{\phantom{l}}_{K_n}/\R_{K_n}^{\gamma^*-1}\mu_{K_n}^{\phantom{lc}})\right)$.}
\end{Sco}

\newpage

\noindent{\bf 3. Discussion du résultat obtenu}\medskip


Une première conséquence de l'expression exacte donnée par le Théorème \ref{RadZ} et son Scolie est de donner immédiatement une généralisation (pour l'exposant $\ell^m$) de l'encadrement proposé par Seo (cf. \cite{Seo1, Seo2}) pour le sous-groupe d'exposant $\ell$ du radical initial des $\Zl$-extensions puis étendu dans \cite{Seo3} au groupe d'exposant $\ell^m$. Il vient, en effet:

\begin{Prop}
On a la double inclusion (où $\pi_0:(\x_n)_{n\in\NN}\mapsto\x_0$ désigne la projection canonique sur le premier terme de la limite projective dans chacun des produits qui suivent):\smallskip

\centerline{$\pi_0 \left( \varprojlim \;(\R^{\phantom{l}}_{K_n}/\R_{K_n}^{\ell^n(\gamma-1)}\R_K^{\ell^m})\right) 
\subset \pi_0 \left( \varprojlim \;(\R^{\phantom{l}}_{K_n}/\R_{K_n}^{\gamma^*-1})\right) 
= \pi_0 \left( \varprojlim \;(\R^{\phantom{l}}_{K_n}/\R_{K_n}^{\gamma^*-1}\R_K^{\ell^m})\right) $.}
\end{Prop}

\Preuve L'égalité de droite est immédiate, puisqu'on a trivialement: $\R_K^{\ell^m}=\R_K^{\gamma^*-1}\subset \R_{K_n}^{\gamma^*-1}$.\smallskip

 Pour établir l'inclusion de gauche, observons d'abord que nous avons de même:

\centerline{ $\R_{K_n}^{\ell^{m+n}}=\R_{K_n}^{\gamma^{*\ell^n}-1}$.}\smallskip

 Partant alors de l'identité évidente $(\gamma-1)-(\gamma^*-1)=\gamma\;\ell^m/(1+\ell^m)\sim\ell^m$, nous obtenons:\smallskip

\centerline{$\R_{K_n}^{(\gamma-1)\ell^n} \subset\; \R_{K_n}^{(\gamma^*-1)\ell^n} \R_{K_n}^{\ell^{m+n}} \subset\; \R_{K_n}^{(\gamma^*-1)\ell^n} \R_{K_n}^{\gamma^{*\ell^n}-1} \subset\;  \R_{K_n}^{\gamma^*-1}$;}\smallskip

\noindent ce qui, joint à l'inclusion $\R_K^{\ell^m}=\R_K^{\gamma^*-1}\subset \R_{K_n}^{\gamma^*-1}$ déjà utilisée, nous montre que les dénominateurs dans la limite projective à gauche sont respectivement contenus dans leurs homologues dans la limite projective au centre. La condition de cohérence requise étant de ce fait plus forte, nous concluons à l'inclusion annoncée.\medskip

\noindent{\bf Nota.} En fait,  Seo (cf. \cite {Seo1, Seo2}) raisonne  sur les quotients $K_n^\times/K_n^{\times\ell^n(\gamma-1)}K^{\times\ell^m}$ pour la limite projective de gauche et $K_n^\times/K_n^{\times(\gamma^*-1)}K^{\times\ell^{m-1}}$ pour celle de droite. Comme ce sont des $\ell$-groupes, on peut les $\ell$-adifier en tensorisant par $\Zl$, ce qui donne respectivement $\R_{K_n}/\R_{K_n}^{\ell^n(\gamma-1)}\R_K^{\ell^m}$ pour le premier, $\R_{K_n}/\R_{K_n}^{\gamma^*-1}\R_K^{\ell^{m-1}}$ pour le second. Le choix malheureux de $\R_K^{\ell^{m-1}}$ dans ce dernier (pour $n>0$) au lieu de $\R_K^{\ell^m}$ a pour conséquence de grossir le dénominateur donc d'affaiblir la condition de cohérence et conduit à une inclusion (potentiellement stricte) au lieu de l'égalité à droite donnée par la Proposition.\medskip

Passant au quotient par les racines de l'unité (qui forment un système projectif pour les applications normes), nous obtenons directement une généralisation des encadrements de Seo:

\begin{Cor}
Le radical initial des $\Zl$-extensions satisfait à l'encadrement:\smallskip

\centerline{$\pi_0 \left( \varprojlim \;\R^{\phantom{l}}_{K_n}/\R_{K_n}^{\ell^n(\gamma-1)}\R_K^{\ell^m}\mu_{K_n}^{\phantom{lc}}\right) 
\subset \pi_0 \left( \varprojlim \;\R^{\phantom{l}}_{K_n}/\R_{K_n}^{\gamma^*-1}\mu_{K_n}^{\phantom{lc}}\right) 
= \pi_0 \left( \varprojlim \;\R^{\phantom{l}}_{K_n}/\R_{K_n}^{\gamma^*-1}\R_K^{\ell^m}\mu_{K_n}^{\phantom{lc}}\right)$.}
\end{Cor}

Oublions maintenant les dénominateurs et considérons la limite projective (toujours pour les applications normes): $\varprojlim \,\R^{\phantom{l}}_{K_n}$. Les suites cohérentes $(x_n)_{n\in\NN}$ qui la constituent sont formées d'idèles principaux $x_n\in\R^{\phantom{l}}_{K_n}$ qui sont en particulier normes dans tous les étages finis $K_{n+k}/K_n$ de la tour $K_\infty/K_n$; et ce sont donc, comme expliqué dans \cite{J28} , des {\em unités logarithmiques}. Il suit:\smallskip

\centerline{$\varprojlim \,\R_{K_n} =\, \varprojlim \,\wE_{K_n}$,}\smallskip

\noindent où $\,\wE_{K_n}$ désigne le $\ell$-groupe des unités logarithmiques du corps $K_n$ (cf. \cite{J28}, \cite{J31} ou \cite{J55}). Les éléments de $\pi_0\left(\varprojlim \,\R^{\phantom{l}}_{K_n}\right)$ sont donc exactement les normes d'unités logarithmiques dans la tour $K_\infty/K$, c'est-à-dire les éléments du sous-groupe $N_{K_\infty\!/K}(\wE_{K_\infty})\underset{\textrm{déf}}{=}\cap_{n\in\NN}\,N_{K_n/K}(\wE_{K_n})$ étudié par Greither (cf. \cite{Grt}). En particulier, il vient ainsi:

\begin{Cor}
Le radical initial des $\Zl$-extensions $\pi_0 \left( \varprojlim \;(\R^{\phantom{l}}_{K_n}/\R_{K_n}^{\gamma^*-1}\mu_{K_n}^{\phantom{lc}})\right)$ contient canoniquement l'image  $\pi_0 \left(\varprojlim \,\R^{\phantom{lc}}_{K_n}\right)\R^{\ell^m}_K/\R_K^{\ell^m}\mu_K^{\phantom{lc}}
=N_{K_\infty\!/K}(\wE^{\phantom{l}}_{K_\infty})\R^{\ell^m}_K/\R_K^{\ell^m}\mu_K^{\phantom{lc}}$ du groupe $N_{K_\infty\!/K}(\wE_{K_\infty})$ des normes cyclotomiques d'unités logarithmiques.
\end{Cor}

Pour l'étude des cas de coïncidence entre $\wE_K$ et son sous-groupe  $N_{K_\infty\!/K}(\wE_{K_\infty})$, on pourra se reporter à \cite{Grt} ou à \cite{J55}.  Pour le lien avec le module de Kuz'min-Tate, voir par exemple \cite{J55}.

\newpage

\noindent{\bf 4. Lien avec le noyau de Tate}\medskip


Les résultats de Tate obtenus par voie cohomologique (cf. \cite{Ta}) montrent que la $\ell$-partie, disons $K_2(K)$ (le premier $\ell$ étant sous-entendu), du groupe universel  pour les symboles sur un corps de nombres $K$ s'obtient simplement par montée et descente dans la $\ell$-tour cyclotomique $K_\infty/K$.\par
Plus précisément, dans le cas qui nous intéresse ici où le corps de base contient les racines $2\ell$-ièmes de l'unité, pour chaque étage fini de la tour $K_\infty/K$, le groupe $K_2(K_n)$ s'écrit canoniquement:\smallskip

\centerline{$K_2(K_n) \simeq (\Bmu^{\phantom{l}}_{\ell^\infty\!} \otimes_{\ZZ} K_\infty^\times)^{\Gamma_{\!n}}/(\Bmu^{\phantom{l}}_{\ell^\infty}\! \otimes_{\ZZ} K_\infty^\times)^{\Gamma_{\!n}}_{_{\mathrm{div}}}$}\smallskip

\noindent comme quotient du sous-groupe des points fixes par $\Gamma_{\!n}=\Gal(K_\infty/K_n)$ du tensorisé $\Bmu^{\phantom{l}}_{\ell^\infty}\! \otimes_{\ZZ} K_\infty^\times$ par son sous-module divisible maximal. Notant alors $\Tl=\varprojlim \Bmu_{\ell^n}$ le module de Tate construit sur les racines de l'unité d'ordre $\ell$-primaire, on peut écrire l'identité précédente sous la forme:\smallskip

\centerline{$K_2(K_n) \simeq (\Tl \otimes_{\Zl} {\mathfrak R}_{K_{\infty\!}})^{\Gamma_{\!n}}/(\Tl \otimes_{\Zl} {\mathfrak R}_{K_{\infty\!}})^{\Gamma_{\!n}}_{_{\mathrm{div}}}$}\smallskip

\noindent (cf. \cite{J18},  I.2, ou \cite{J23}, Th. 1.4); soit encore, en posant $\bar\gamma^*= \gamma^{-1}/(1+\ell^m)$ et $\bar\Gamma^*_{\!n}=\bar\gamma^{*\,\ell^n\Zl}$:\smallskip

\centerline{$K_2(K_n) \simeq \Tl \otimes_{\Zl} \left({\mathfrak R}_{K_\infty}^{\,\bar\Gamma\*_n}/({\mathfrak R}_{K_{\infty\!}}^{\,\bar\Gamma\*_n})_{_{\mathrm{div}}}\right)$}.\smallskip

\noindent Le sous-groupe divisible $({\mathfrak R}_{K_{\infty\!}}^{\,\bar\Gamma\*_n})_{_{\mathrm{div}}}$ est ainsi le noyau universel de Tate. Par substitution de $\bar\gamma=\gamma^{-1}$ à $\gamma$, il se présente formellement tout comme le radical des $\Zl$-extensions, de sorte que nous pouvons lui appliquer directement {\em mutatis mutandis} les mêmes arguments pour obtenir:

\begin{Th}
Le sous-module d'exposant $\ell^m$du noyau universel de Tate attaché au corps $K$, i.e. le sous-groupe du radical universel ${\mathfrak R}_K$ défini par \smallskip

\centerline{$_{\ell^m\!}{\mathfrak N}_K = \{\ell^{-m}\otimes x \in {\mathfrak R}_K \; | \;  \{\zeta_{\ell^n},x\}=1\;\mathrm{dans}\;K_2(K)\}$,}\smallskip 

\noindent est formé des termes initiaux $\x_0$ de la limite projective pour les applications normes $\varprojlim \,{\mathfrak  R}_{K_n}^{\,\bar\Gamma^*}$.\par

Il s'identifie donc au terme initial de la limite projective des quotients des groupes d'idèles principaux invariants par $\bar\gamma^*-1$ attachés aux étages finis de la tour cyclotomique $K_\infty/K$:
\smallskip

\centerline{$_{\ell^m\!}{\mathfrak N}_K = \pi_0 \left( \varprojlim \;(\R^{\phantom{l}}_{K_n}/\R_{K_n}^{\bar\gamma^*-1}\mu_{K_n}^{\phantom{lc}})\right)$.}
\end{Th}

\begin{Cor}
Le noyau universel de Tate  $_{\ell^m\!}{\mathfrak N}_K$ vérifie en particulier l'encadrement:\smallskip

\centerline{$\pi_0 \left( \varprojlim \;\R^{\phantom{l}}_{K_n}/\R_{K_n}^{\ell^n(\gamma-1)}\R_K^{\ell^m}\mu_{K_n}^{\phantom{lc}}\right) 
\subset \pi_0 \left( \varprojlim \;\R^{\phantom{l}}_{K_n}/\R_{K_n}^{\bar\gamma^*-1}\mu_{K_n}^{\phantom{lc}}\right) 
= \pi_0 \left( \varprojlim \;\R^{\phantom{l}}_{K_n}/\R_{K_n}^{\bar\gamma^*-1}\R_K^{\ell^m}\mu_{K_n}^{\phantom{lc}}\right)$.}
\end{Cor}

On notera que le groupe à gauche est identique à celui obtenu précédemment. En particulier:

\begin{Cor}
Le noyau $_{\ell^m\!}{\mathfrak N}_K$ contient canoniquement l'image  $\pi_0 \left(\varprojlim \,\R^{\phantom{l}}_{K_n}\right)\R^{\ell^m}_K/\R_K^{\ell^m}\mu_K^{\phantom{lc}}=N_{K_\infty\!/K}(\wE_{K_\infty})\R^{\ell^m}_K/\R_K^{\ell^m}\mu_K^{\phantom{lc}}$ du groupe $N_{K_\infty\!/K}(\wE_{K_\infty})$ des normes  d'unités logarithmiques.
\end{Cor}

Pour conclure, précisons quelques points:\smallskip

(i)  En définissant pour tout $i\in\ZZ$ le $i$-ième tordu à la Tate du générateur topologique $\gamma$ par $\gamma^{(i)}=\gamma/(1+\ell^m)^i$ et en introduisant le groupe $\Gamma^{(i)}=\gamma^{(i)\,\Zl}$, on obtient facilement des résultats  en tous points analogues à ceux présentés plus haut pour le radical initial des $\Zl$-extensions (qui correspond au choix $i=1$) et pour le noyau des symboles universels (qui correspond, lui, au choix $i=-1$). Les résultats et conjectures standard de Schneider liées à ces divers noyaux sont présentées dans \cite{JM} en lien avec l'arithmétique des classes logarithmiques.\smallskip

(ii) Il est toujours possible, du moins pour $\ell$ impair, de s'affranchir de l'hypothèse $\Bmu^{\phantom{l}}_\ell \subset K$: il suffit, en effet, d'introduire le corps $K'=K[\zeta_\ell]$, de raisonner sur $K'$, puis de redescendre les résultats sur $K$ en utilisant les idempotents de l'algèbre semi-locale $\Zl[\Delta]$ associés aux caractères $\ell$-adiques du groupe $\Delta=\Gal(K'/K)$. Il convient de prendre garde, en revanche, au fait que la torsion par les racines de l'unité introduit un décalage des caractères et que la descente ne se fait donc pas toujours via le caractère unité (cf. e.g. \cite{J28} ou \cite{JM} pour plus de détails).\smallskip

(iii) Il peut y avoir égalité ou pas  entre le noyau de Tate $_{\ell^m\!}{\mathfrak N}_K$  et le radical initial $_{\ell^m\!}{\mathfrak Z}_K$, comme le montrent, par exemple, les tables de Thomas (cf. \cite{Th}). L'égalité a lieu en particulier (mais sans que cette condition soit nécessaire) lorsque le $\ell$-groupe des classes logarithmiques est trivial, auquel cas les deux groupes coïncident avec le radical hilbertien étudié dans \cite{J23} et {\cite{J24}.

\newpage

\def\refname{\normalsize{\sc  Références}}

\bigskip\noindent
{\small
\begin{tabular}{l}
{Jean-François {\sc Jaulent}}\\
Institut de Mathématiques de Bordeaux \\
Université de {\sc Bordeaux} \& CNRS \\
351, cours de la libération\\
F-33405 {\sc Talence} Cedex\\
courriel : Jean-Francois.Jaulent@math.u-bordeaux1.fr 
\end{tabular}
}

 \end{document}